%% file: conference_101719.tex
\documentclass[conference]{IEEEtran}
\IEEEoverridecommandlockouts
\usepackage{cite}
\usepackage{amsmath,amssymb,amsfonts}
\usepackage{algorithmic}
\usepackage{graphicx}
\usepackage{textcomp}
\usepackage{xcolor}
\usepackage{comment}
\usepackage{todonotes}
\usepackage{mathtools}
\usepackage{cleveref}
\usepackage{siunitx}

\input{Macros.tex}

\input{Packages.tex}

\def\BibTeX{{\rm B\kern-.05em{\sc i\kern-.025em b}\kern-.08em
    T\kern-.1667em\lower.7ex\hbox{E}\kern-.125emX}}
\begin{document}

\title{On the Fast Direct Solution of a Preconditioned Electromagnetic Integral Equation\\
\thanks{This work has been funded by the European Research Council (ERC) under the European Union’s Horizon 2020 research and innovation program (ERC project 321, grant No.724846) and by the ANR Labex CominLabs under the project ``CYCLE''.}
}

\author{Davide~Consoli\textsuperscript{1}, Clément~Henry\textsuperscript{2}, Alexandre~Dély\textsuperscript{1}, Lyes~Rahmouni\textsuperscript{1}, John~Erik~Ortiz~Guzman\textsuperscript{3},\\Tiffany~L.~Chhim\textsuperscript{1}, Simon~B.~Adrian\textsuperscript{4}, Adrien~Merlini\textsuperscript{2}, and Francesco~P.~Andriulli\textsuperscript{1} \\
\textsuperscript{1}Politecnico di Torino, Turin, Italy\\
\textsuperscript{2}IMT Atlantique, Brest, France\\
\textsuperscript{3}Universidad de Nariño, Pasto, Colombia\\
\textsuperscript{4}Universität Rostock, Rostock, Germany}

\maketitle

\begin{abstract}
This work presents a fast direct solver strategy for electromagnetic integral equations in the high-frequency regime. 
The new scheme relies on a suitably preconditioned combined field formulation and results in a single skeleton form plus identity equation. This is obtained after a regularization of the elliptic spectrum  through the extraction of a suitably chosen equivalent circulant problem. The inverse of the system matrix is then obtained by leveraging the Woodbury matrix identity, the low-rank representation of the extracted part of the operator, and fast circulant algebra yielding a scheme with a favorable complexity and suitable for the solution of multiple right-hand sides. Theoretical considerations are accompanied by numerical results both of which are confirming and showing the practical relevance of the newly developed scheme.  
\end{abstract}


\section{Introduction}

Integral equation strategies are effective for solving scattering and radiation problems over a wide range of frequencies since they do not require absorbing boundary conditions, they are free from numerical dispersion, and they only require the discretization of the surfaces of the scatterers and the radiators.
Unfortunately, matrices arising from their discretization are dense and have dimensions that grow for increasing frequencies. Naive solutions would scale cubically with the number of unknowns and easily jeopardize the advantages of integral equation approaches. Fast matrix-vector multiplication algorithms \cite{chew_fast_2001} can substantially decrease this computational burden, but they require iterative solutions, which renders them less attractive when many right-hand sides (RHSs) must be considered. An alternative is the use of fast direct solvers that are building, in reduced complexity, directly the inverse of the system matrix and are efficient for solving problems with multiple RHSs. There are several effective strategies for direct solutions that often rely on hierarchical decompositions both for problems in the low- and in the high-frequency regimes (see \cite{adams_modular_2008,guo_butterfly-based_2017,sharshevsky_direct_2020} and references therein).

In this work, we leverage a different approach. We first precondition the combined field integral equation (CFIE) in a suitable way \cite{andriulli_high_2015}. We then exploit the particular structure of the preconditioned spectrum to separate its principal part contribution and compress its remainder. This is made practical by extracting an equivalent circulant
problem, which has the advantage of automatically extracting the principal part contribution and the remainder is compressed in favorable complexity. This results in a non-hierarchical scheme in which a single skeleton form is required and the solution for several RHSs can be obtained efficiently. Theoretical considerations are matched with numerical results confirming and showing the practical relevance of the new scheme. 

\section{Notation and Background}

Consider the scattering of a TE time-harmonic incident electromagnetic field ($\veg E^\text{inc}$, $\veg H^\text{inc}$) of angular frequency $\omega$ impinging on a perfectly electrically conducting (PEC) body, modeled by a 2-dimensional convex contour $\gamma$ that resides in a medium of permittivity $\epsilon$, permeability $\mu$, impedance $\eta = (\mu/\epsilon)^{1/2}$, and corresponding wavenumber $k=\omega \sqrt{\mu \epsilon}$. The combined field integral equation (CFIE) that relates the induced tangential  current density $\veg j$ on $\gamma$ to the incident electric and magnetic fields reads
\begin{equation}\label{eq:cont}
    \frac{1}{\im k} \op N^k (j_t) + \left(\frac{\op I}{2} - \op D^k\right)(j_t) = - \frac{1}{\eta} E_t^\text{inc} - H_z^\text{inc}
\end{equation}
with 
\begin{align}
  (\op N^k j_t) (r) &\coloneqq -\frac{\partial}{\partial n} \int_{\gamma} \frac{\partial}{\partial n'} g_k(\vv r, \vv r') j_t(r') \dd r'\,,\\
  (\op D^k  j_t) (r) &\coloneqq  \int_{\gamma} \frac{\partial}{\partial n'} g_k(\vv r, \vv r') j_t(r') \dd r'\,,
\end{align}
and $g_k(\vv r, \vv r') = \im/4 H^{(1)}_0(k \|\vv r - \vv r'\|)$. 
 After boundary element discretization with piecewise linear functions $\lambda_i(\vv r)$, \eqref{eq:cont} is discretized as
\begin{equation}
     \frac{1}{\im k} \mat N^k \mat j + \left(\frac{\mat G}{2} - \mat D^k\right) \mat j= -\frac{1}{\eta}\mat e - \mat h
\end{equation}
with $\left[\mat e\right]_{i} = \int_\gamma \lambda_i(\vv r) E^\text{inc}_t(\vv r) \dd \vv r$,$\left[\mat h\right]_{i} = \int_\gamma \lambda_i(\vv r) H^\text{inc}_z(\vv r) \dd \vv r$, $[\mat N^k]_{ij} = \left<\lambda_i, \op N^k \lambda_j\right>$, $[\mat D^k]_{ij} = \left<\lambda_i, \op D^k \lambda_j\right>$, $[\mat G]_{ij} = \left<\lambda_i, \lambda_j\right>$, 
with
$\left<f,g\right> = \int_\gamma f(\vv r) g(\vv r) \dd \vv r$.


\section{A Suitable Integral Equation and its Spectral Analysis}

\begin{figure}
\centering
\includegraphics[width=1\linewidth]{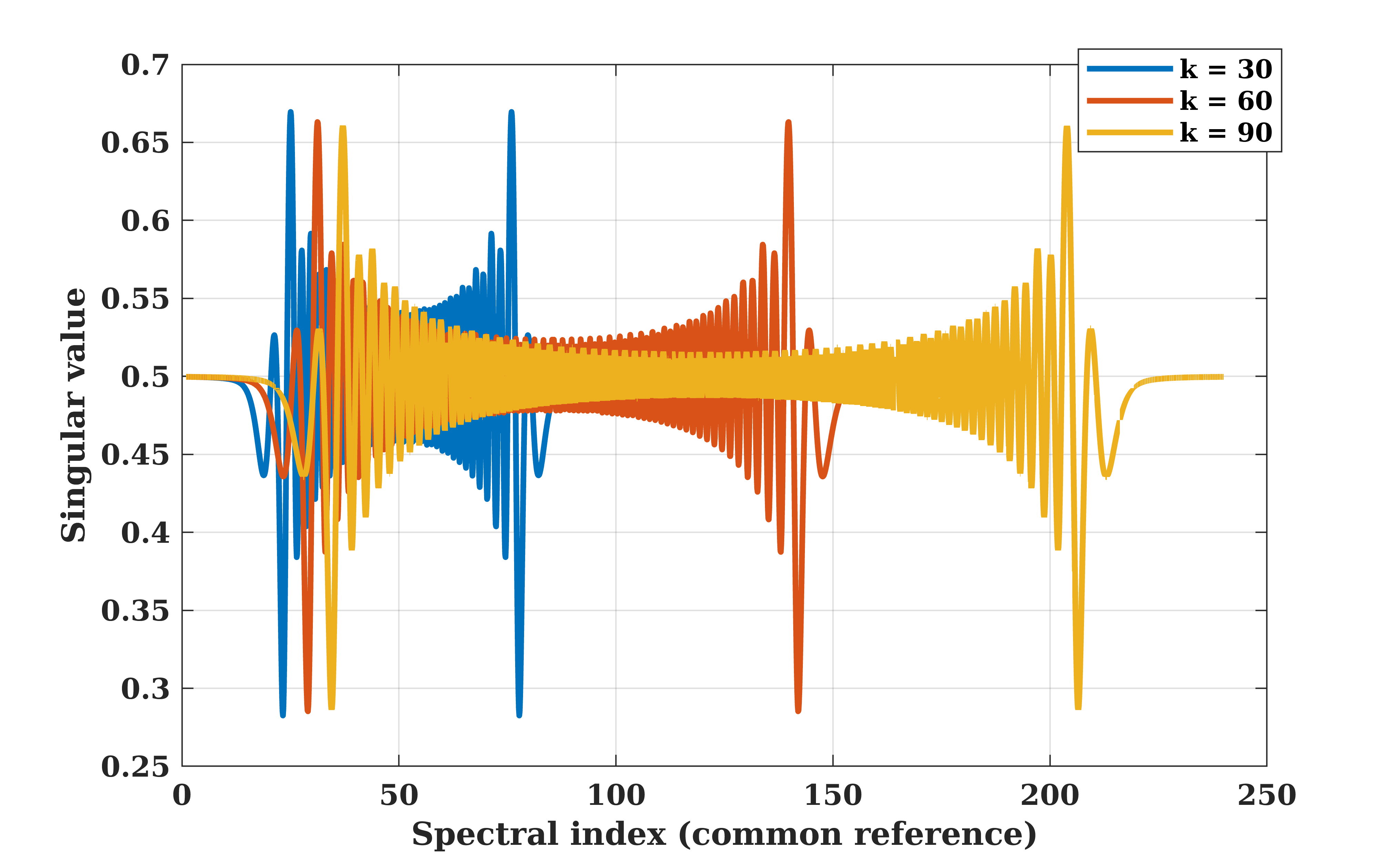}
\caption{Spectrum of the equation of this work (continuous counterpart of \eqref{eq:equationprecond}) obtained for a cylinder.}
\label{fig:spectrumthiswork}
\end{figure}

\begin{figure}
\centering
\includegraphics[width=1\columnwidth]{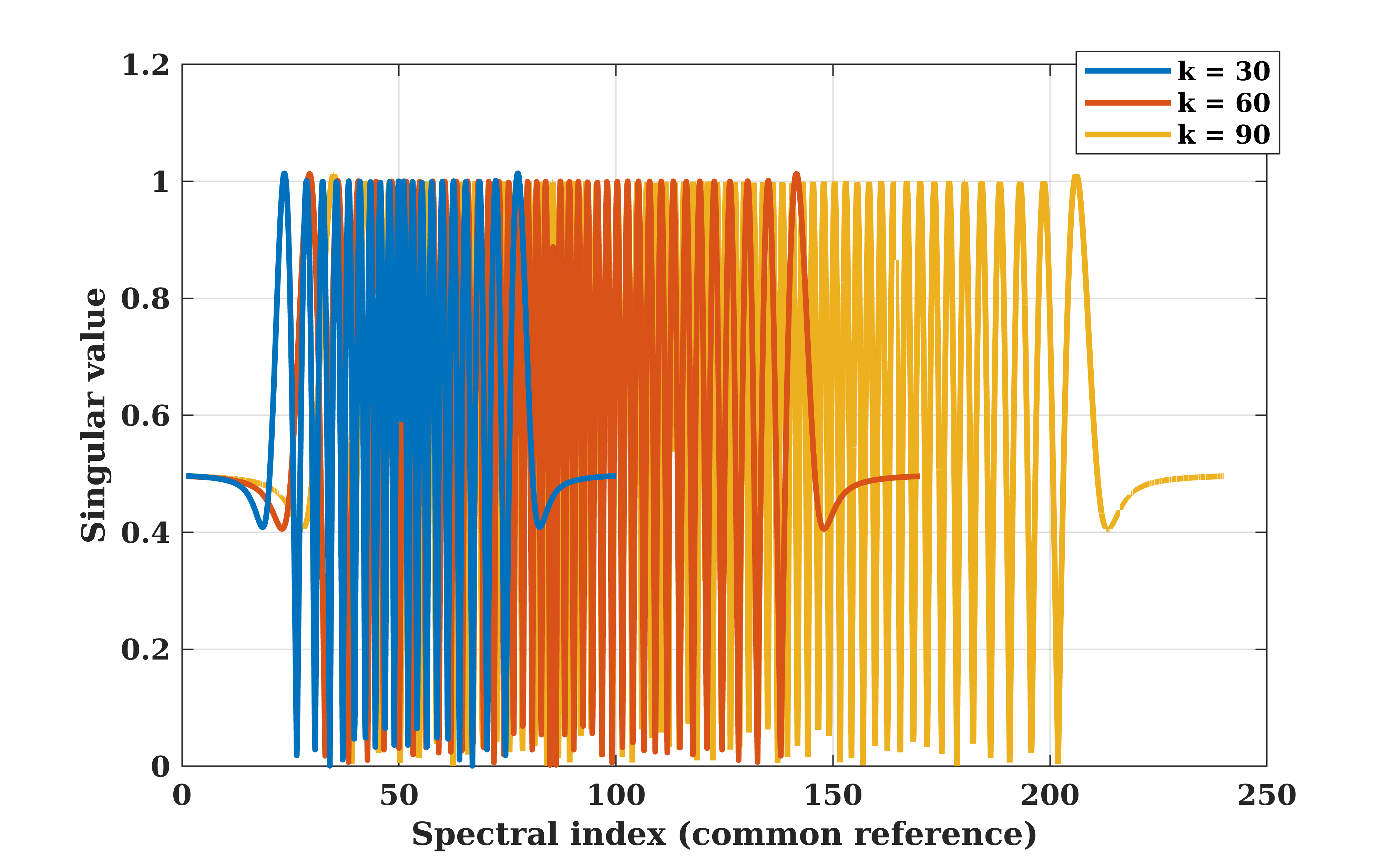}
\caption{Continuous spectrum of the 2D TE-MFIE obtained for a cylinder.}
\label{fig:MFIE}
\end{figure}

A well conditioned equation can be obtained from \eqref{eq:cont} by leveraging a modified version of the preconditioner proposed in \cite{andriulli_high_2015}
\begin{multline}
    \label{eq:equationprecond}
    \left(\mat S^{\tilde{k}}\mat G^{-1}\mat N^{k}+\left(\frac{\mat G}{2}+\mat D^{\tilde{k}}\right)\mat G^{-1}\left(\frac{\mat G}{2}-\mat D^k\right)\right)\mat j \\= -\frac{\im k }{\eta}\mat S^{\tilde{k}}\mat G^{-1}\mat e - \left(\frac{\mat G}{2}+\mat D^{\tilde{k}}\right)\mat G^{-1} \mat h\,,
\end{multline}
where $\tilde k \coloneqq k + 0.4\im k^{1/3} a^{-2/3}$ following \cite{darbas_generalized_2006,boubendir_well-conditioned_2014} with $a$ a suitable average of the radii of curvature along $\gamma$.
When comparing the spectrum of the above equation (Fig.~\ref{fig:spectrumthiswork}) with that of the standard formulations such as the EFIE, the MFIE (Fig.~\ref{fig:MFIE}), or even the Calderón-preconditioned CFIE, one can see that the spectral content of \eqref{eq:equationprecond} clusters, and is maximal, around the surface resonant point attained when the spatial frequency equals $k a$. Note that for ease of reading the singular values are ordered by the Fourier content of their singular vectors. This suggests that selecting the singular vectors corresponding to this region of maximal spectral strength and filtering out the deviation from the halved identity of the others would lead to a compression of the electromagnetic operator in the high frequency regime---in which the discretization parameter is kept at a fixed proportion of the wavelength. In particular, this would correspond to searching the solution in the space of corrections to the RHS in \eqref{eq:equationprecond} of the form $s\mapsto M(s)\e^{\im k s}$, with $M(s)$ slowly varying.

\section{Handling the Error on the Elliptic Spectrum}

\begin{figure}
\centering
\includegraphics[width=1\linewidth]{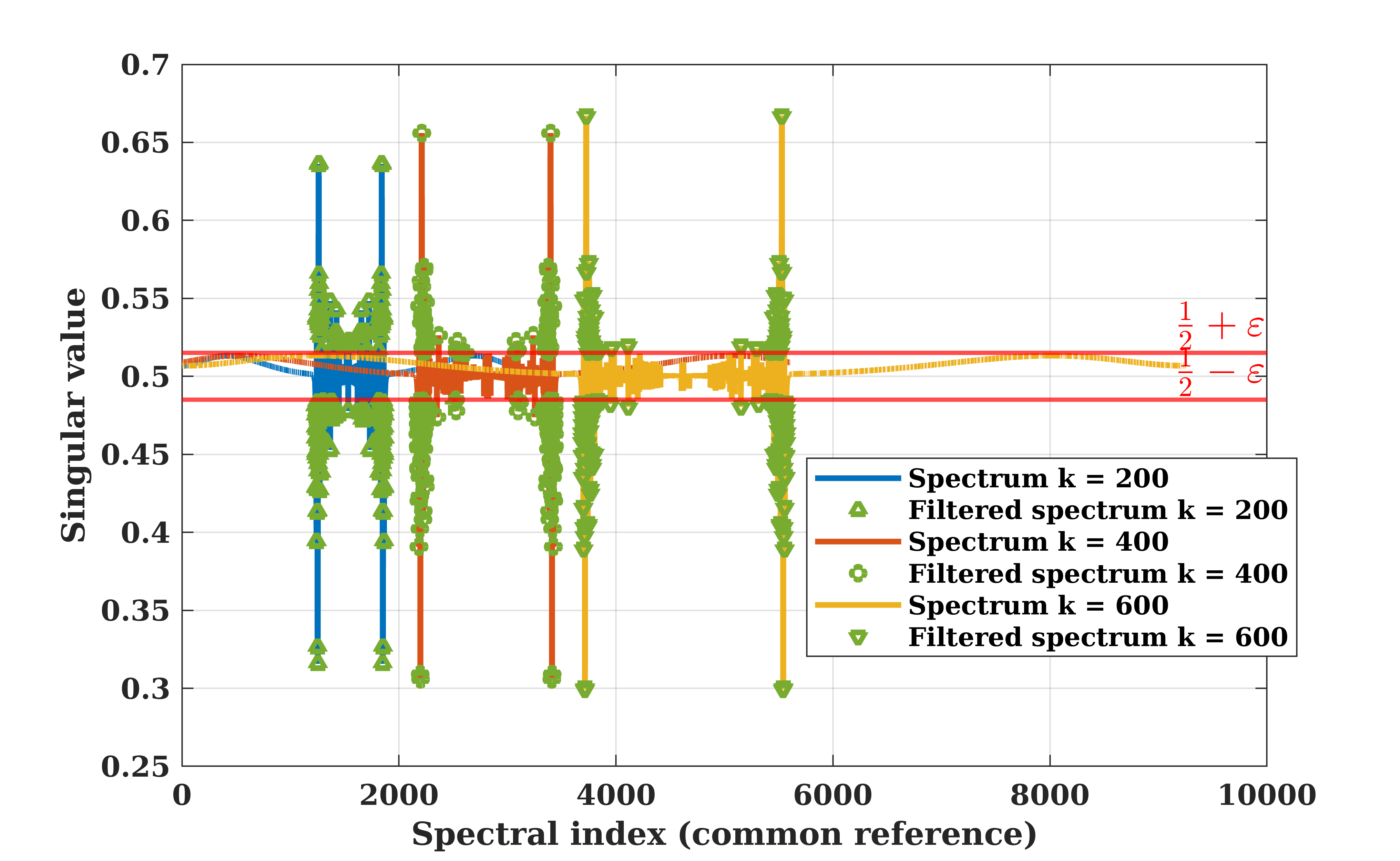}
\caption{Spectrum of the matrix $\mat C$ for increasing frequency, with a compression parameter $\varepsilon = 0.015$. The green markers indicate which singular values/singular vectors are selected by the compression scheme.}
\label{fig:compressionNonextracted}
\end{figure}

Although  eigenvector extraction (after proper shifting and symmetrization) is an option for obtaining a direct solver out of \eqref{eq:equationprecond}, at higher precision this may become less practical and skeleton form obtaining algorithms can be used instead. However, the discretization error in the elliptic part of the spectrum can interfere with this approach. The discretization will cause a constant ($\mathcal O(1)$) error in the elliptic spectrum for first kind preconditioned operators (see Fig.~\ref{fig:compressionNonextracted}).
There are several strategies to  overcome this problem. 
Here we opted for subtracting an equivalent circular problem from the operators defined on $\gamma$.
This will, in practice, extract the second kind part of the original operators, leaving only a compact operator for which the relative error (with respect to the identity) will be decreasing.

\section{The Direct Solver}

To exploit these considerations to the fullest, we will study the deviation of the electromagnetic operators defined on $\gamma$ from their equivalents defined on the circle $\gamma_c$ that has the same perimeter as $\gamma$. These auxiliary matrices, denoted by the subscript ${}_c$, are obtained with the same uniform discretization in the curvilinear abscissa as the main matrices defined on $\gamma$---both of which are discretized with piecewise linear basis functions in the curvilinear abscissa.
Note that under these conditions we have $\mat G = \mat G_c$, among other critical properties.
Then by defining
\begin{equation}\label{eq:LHSequationprecond}
    \mat C\coloneqq\left(\left(\frac{\mat G}{2}+\mat D^{\tilde{k}}\right)\mat G^{-1}\left(\frac{\mat G}{2}-\mat D^k\right)+\mat S^{\tilde{k}}\mat G^{-1}\mat N^{k}\right)
\end{equation}
and its circle counterpart $\mat C_c$, we compute the skeleton form of the matrix
\begin{equation}
    \mat C-\mat C_c\eqqcolon\mat U \mat V^\T\,, \label{eq:skeleton}
\end{equation}
which can be obtained efficiently based on the expected low rank of the difference and on the fact that all operations involving the matrices and operators on $\gamma_c$ are done rapidly via fast Fourier transform (FFT).
Finally, the solution of the original problem can be obtained efficiently, for any number of RHSs, by inverting
\begin{equation}
    \mat C = \mat C_c\left( \mat I + \left(\mat C_c\right)^{-1}\mat U \mat V^T\right)\, ,
\end{equation}
where $\mat I$ is the identity matrix, using the Woodbury matrix identity \cite{henderson_deriving_1981}
\begin{equation}
    \mat C^{-1}= \mat C_c^{-1}- \mat C_c^{-1} \mat U \left( \mat I +\mat V^\T \mat C_c^{-1}\mat U \right)^{-1} \mat V^\T\mat C_c^{-1}
\end{equation}
that can be computed efficiently given the low-rank (i.e., the inner dimension) nature of the skeleton form and circulant algebra.

\section{Numerical Results}

\begin{figure}
\centering
\includegraphics[width=1\linewidth]{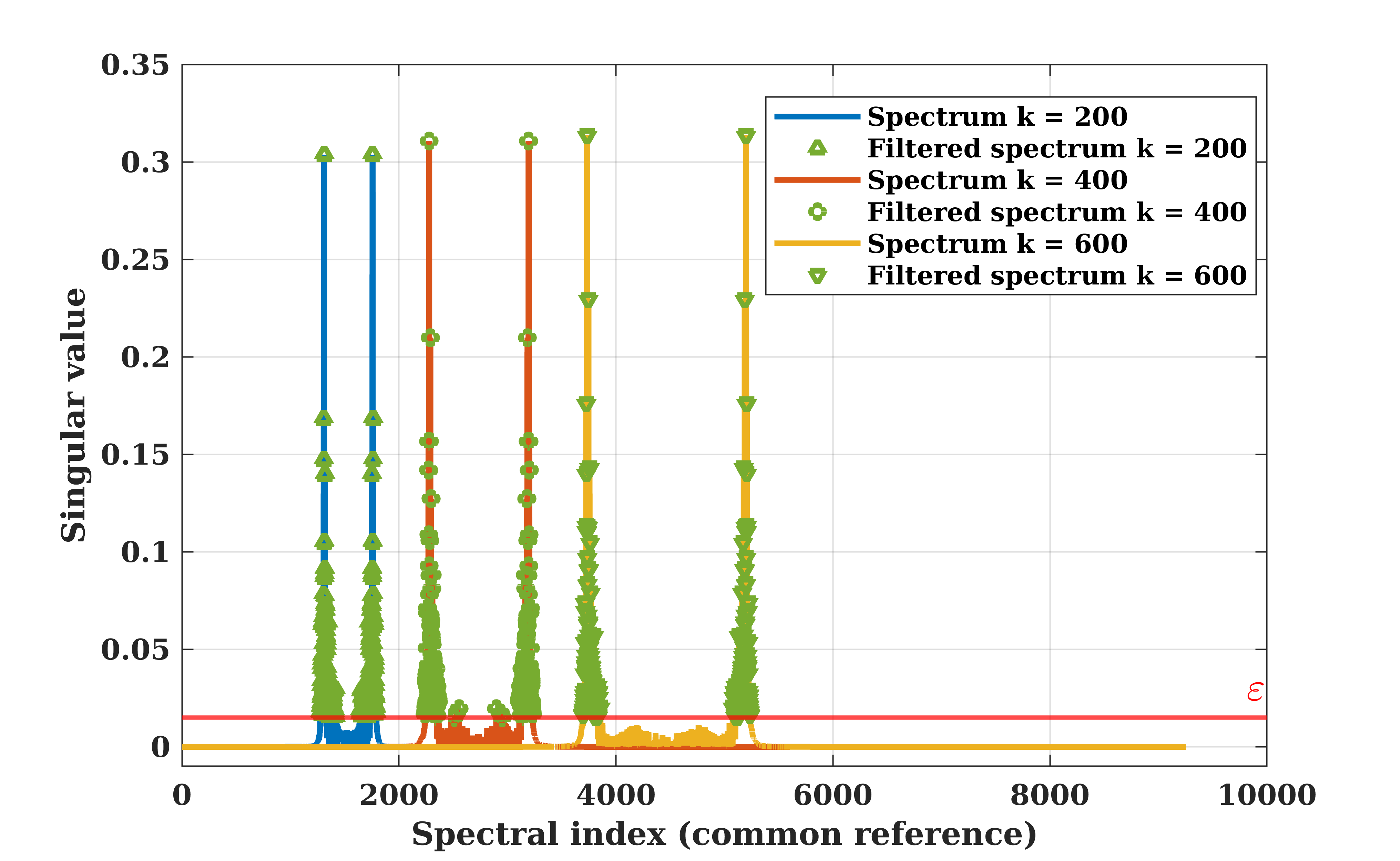}
\caption{Spectrum of the $\mat C - \mat C_c$ for increasing frequencies.}
\label{fig:compressionextracted}
\end{figure}

\begin{figure}
\centering
\includegraphics[width=1\linewidth]{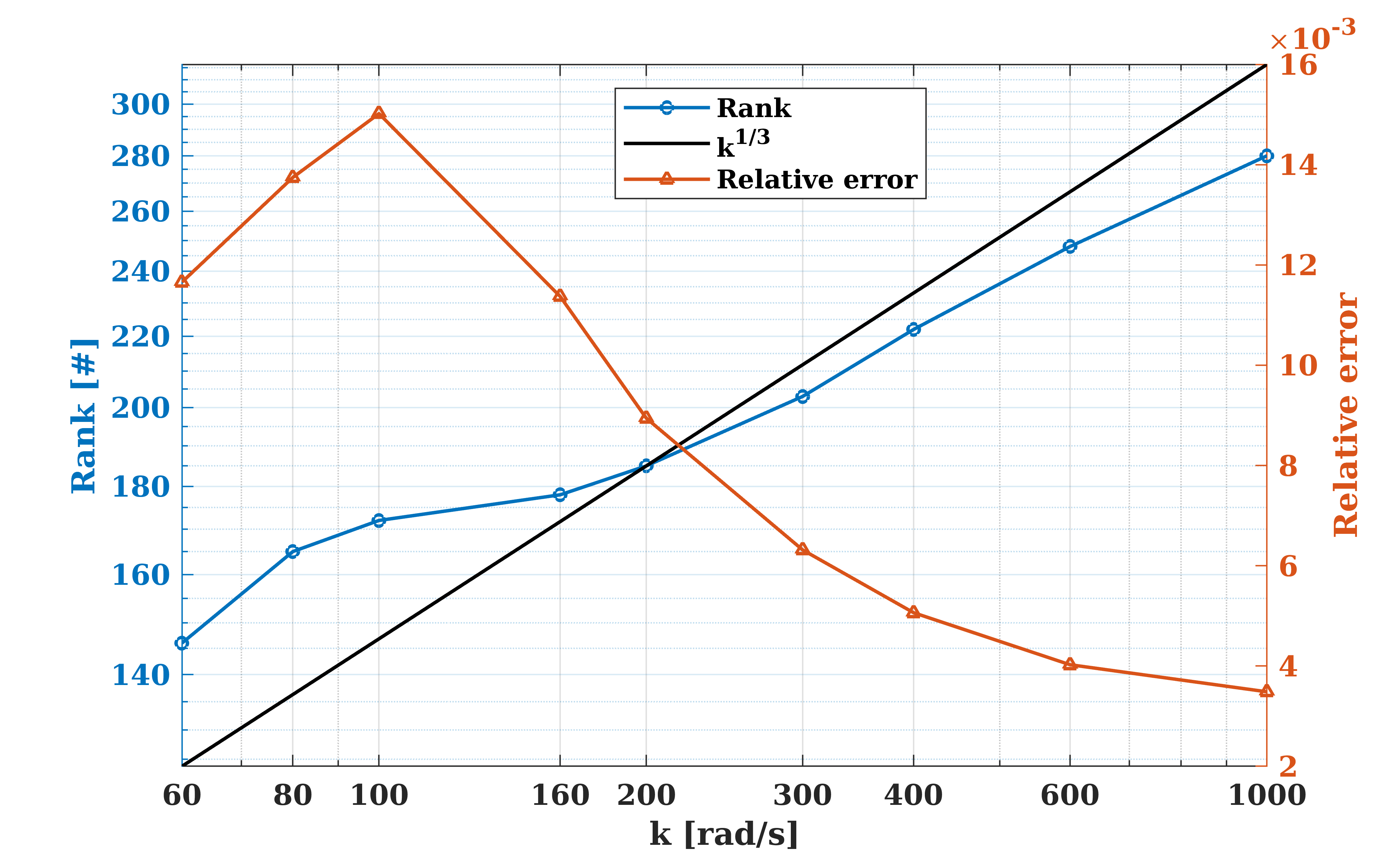}
\caption{Accuracy of the solution and rank of the compressed matrix as a function of $k$, with a compression parameter $\varepsilon = 0.015$.}
\label{fig:precision}
\end{figure}

The direct solution strategy of this work has been tested on  an ellipse with semi-major axis \SI{2}{\meter} and semi-minor axis \SI{1}{\meter}. A close inspection of the spectrum of $\mat C$ (Fig.~\ref{fig:compressionNonextracted}) reveals that the spectral corruption caused by the discretization error poses a real challenge for compression purposes, since they would artificially increase the rank of the compression if the looked-for accuracy was too high (the limit case is indicated by the red lines in the figure). On the contrary, the elliptic spectrum of $(\mat C-\mat C_c)$ (Fig.~\ref{fig:compressionextracted}) does not remain constant, but it decreases for the high spatial frequencies, which makes it an ideal candidate for applying the proposed compression scheme. Fig.~\ref{fig:precision} shows the behavior of the inner rank of the skeleton in \eqref{eq:skeleton} that shows a growth not larger than $k^{\frac{1}{3}}$ (which is a bound that can be theoretically predicted for the circular case). This provides in this case an upper bound for the computational complexity not larger than $O(N^{\frac{4}{3}})$ where in practice lower complexities can be obtained given that for a constant truncation error on the matrix the accuracy of the solution increases with frequency, as can be seen in Fig.~\ref{fig:precision}. Thus, the effectiveness of the scheme is ensured in the high-frequency regime.

\bibliographystyle{IEEEtran}
\bibliography{bibliography}

\end{document}

%% file: Macros.tex
\newcommand{\mr}{\mathrm}                   

\newcommand{\vt}[1]{\boldsymbol{#1}} 

\newcommand{\veg}[1]{\bm{#1}}               
\newcommand{\mat}[1]{\mathbf{#1}}       
\newcommand{\op}[1]{\mathcal{#1}}           

\newcommand{\vv}[1]{\vt{#1}}

\newcommand{\dd}{\mathrm{d}}  

\newcommand{\im}{\mathrm{i}}  

\newcommand{\e}{\mathrm{e}}

\newcommand{\T}{\mr{T}}

%% file: Packages.tex
\usepackage[utf8]{inputenc}

\usepackage[T1]{fontenc}
\usepackage{subcaption}
\usepackage{float}
\usepackage{tikz}
\usepackage{pgfplots}
\usepackage{pgfplotstable}
\usepackage{siunitx}
\usepackage{cleveref}
\pgfplotsset{compat=newest}

\usepackage{amsthm}
\DeclareMathAlphabet{\mathbfsf}{\encodingdefault}{\sfdefault}{bx}{n}
\usepackage{bm}
\usepackage[OMLmathsfit]{isomath}


%% file: conference_101719.bbl
\begin{thebibliography}{1}
\providecommand{\url}[1]{#1}
\csname url@samestyle\endcsname
\providecommand{\newblock}{\relax}
\providecommand{\bibinfo}[2]{#2}
\providecommand{\BIBentrySTDinterwordspacing}{\spaceskip=0pt\relax}
\providecommand{\BIBentryALTinterwordstretchfactor}{4}
\providecommand{\BIBentryALTinterwordspacing}{\spaceskip=\fontdimen2\font plus
\BIBentryALTinterwordstretchfactor\fontdimen3\font minus
  \fontdimen4\font\relax}
\providecommand{\BIBforeignlanguage}[2]{{%
\expandafter\ifx\csname l@#1\endcsname\relax
\typeout{** WARNING: IEEEtran.bst: No hyphenation pattern has been}%
\typeout{** loaded for the language `#1'. Using the pattern for}%
\typeout{** the default language instead.}%
\else
\language=\csname l@#1\endcsname
\fi
#2}}
\providecommand{\BIBdecl}{\relax}
\BIBdecl

\bibitem{chew_fast_2001}
W.~C. Chew, J.-M. Jin, E.~Michielssen, and J.~M. Song, Eds., \emph{Fast and
  {{Efficient Algorithms}} in {{Computational Electromagnetics}}}.\hskip 1em
  plus 0.5em minus 0.4em\relax {Artech House}, 2001.

\bibitem{adams_modular_2008}
R.~J. Adams, Y.~Xu, X.~Xu, J.-s. Choi, S.~D. Gedney, and F.~X. Canning,
  ``Modular {{Fast Direct Electromagnetic Analysis Using Local-Global Solution
  Modes}},'' \emph{IEEE Transactions on Antennas and Propagation}, vol.~56,
  no.~8, pp. 2427--2441, Aug. 2008.

\bibitem{guo_butterfly-based_2017}
H.~Guo, Y.~Liu, J.~Hu, and E.~Michielssen, ``A {{Butterfly-Based Direct
  Integral-Equation Solver Using Hierarchical LU Factorization}} for
  {{Analyzing Scattering From Electrically Large Conducting Objects}},''
  \emph{IEEE Transactions on Antennas and Propagation}, vol.~65, no.~9, pp.
  4742--4750, Sep. 2017.

\bibitem{sharshevsky_direct_2020}
A.~Sharshevsky, Y.~Brick, and A.~Boag, ``Direct {{Solution}} of {{Scattering
  Problems Using Generalized Source Integral Equations}},'' \emph{IEEE
  Transactions on Antennas and Propagation}, vol.~68, no.~7, pp. 5512--5523,
  Jul. 2020.

\bibitem{andriulli_high_2015}
F.~P. Andriulli, I.~Bogaert, and K.~Cools, ``On the high frequency behavior and
  stabilization of a preconditioned and resonance-free formulation,'' in
  \emph{2015 {{International Conference}} on {{Electromagnetics}} in {{Advanced
  Applications}} ({{ICEAA}})}, Sep. 2015, pp. 1321--1324.

\bibitem{darbas_generalized_2006}
M.~Darbas, ``Generalized combined field integral equations for the iterative
  solution of the three-dimensional {{Maxwell}} equations,'' \emph{Applied
  Mathematics Letters}, vol.~19, no.~8, pp. 834--839, Aug. 2006.

\bibitem{boubendir_well-conditioned_2014}
Y.~Boubendir and C.~Turc, ``Well-conditioned boundary integral equation
  formulations for the solution of high-frequency electromagnetic scattering
  problems,'' \emph{Computers \& Mathematics with Applications}, vol.~67,
  no.~10, pp. 1772--1805, Jun. 2014.

\bibitem{henderson_deriving_1981}
H.~V. Henderson and S.~R. Searle, ``On {{Deriving}} the {{Inverse}} of a
  {{Sum}} of {{Matrices}},'' \emph{SIAM Review}, vol.~23, no.~1, pp. 53--60,
  Jan. 1981.

\end{thebibliography}
